\documentclass[a4paper,11pt]{article}
\usepackage[nointlimits]{amsmath}
\usepackage{amssymb,amsmath,amsthm}
\usepackage{amsfonts}
\usepackage{latexsym}
\usepackage{graphicx}
\usepackage{color}
\usepackage{layout}
\usepackage[english]{babel}
\numberwithin{equation}{section} 

\begin{document}

\title{\bf An elementary proof of \\
a Rodriguez-Villegas supercongruence}
\author{{\sc Roberto Tauraso}\\
Dipartimento di Matematica\\
Universit\`a di Roma ``Tor Vergata'', Italy\\
{\tt tauraso@mat.uniroma2.it}\\
{\tt http://www.mat.uniroma2.it/$\sim$tauraso}
}

\date{}
\maketitle
\begin{abstract}
\noindent We give a short proof of the following known congruence: for every odd prime $p$
$$\sum_{k=0}^{p-1}{2k\choose k}^2 16^{-k}\equiv (-1)^{{p-1\over 2}}\pmod{p^2}.$$
Moreover, we provide some new results connected with the above congruence.
\end{abstract}

\makeatletter{\renewcommand*{\@makefnmark}{}
\makeatother}

\section{Introduction} 
The congruence mentioned in the abstract had been conjectured by Rodriguez-Villegas in \cite{RV:03} and it is
connected with some properties of the hypergeometric families of Calabi-Yau manifolds. 
In \cite{Mo:03}, Mortenson proved the congruence by using the theory of Gaussian hypergeometric series, the properties of the $p$-adic
$\Gamma$-function, and a strange combinatorial identity. 
In this short note we would like to present a new proof of the
congruence which is elementary and uses a different combinatorial identity.
The following theorem establishes a more general result, and the conjecture follows immediately by taking $\alpha=0$ and $\beta=-1$
and by noting that $p^2$ divides ${2k\choose k}^2$ for ${p-1\over 2}<k<p-1$.

\vspace*{2mm}

\noindent {\bf Theorem} {\sl Let $p$ be an odd prime and let $\alpha,\beta\in\mathbb{Z}_p$, then
$$\sum_{k=0}^{{p-1\over 2}}
{2k\choose k}^2 {(\alpha-\beta)^{{p-1\over 2}-k}(-\beta)^k\over 16^k}\equiv
\sum_{k=0}^{{p-1\over 2}}{{p-1\over 2}\choose k}^2\alpha^{{p-1\over 2}-k}\beta^k\pmod{p^2}.$$}

\vspace*{2mm}

\noindent {\bf Remarks.} If we denote by $T(\alpha,\beta)$ the LHS of the above congruence then it is easy
to realize (see the RHS) that this function has some interesting symmetries
$$T(\alpha,\beta)\equiv T(\beta,\alpha)\equiv (-1)^{{p-1\over 2}}T(-\alpha,-\beta)\pmod{p^2}.$$
Moreover, by setting in a proper way the parameters $\alpha$ and $\beta$ we can show several results:

\noindent (i) By letting $\alpha=1-t$ and $\beta=-t$ then by the above symmetry we have that 
$$\sum_{k=0}^{p-1}{2k\choose k}^2 \left({t\over 16}\right)^k\equiv 
(-1)^{{p-1\over 2}}\sum_{k=0}^{p-1} {2k\choose k}^2 \left({1-t\over 16}\right)^k\pmod p^2.$$
Hence, since the Fibonacci and the Lucas sequences are respectively given by
$F_k=(t^k-(1-t)^k)/\sqrt{5}$ and $L_k=(t^k+(1-t)^k)$
for $t=(1+\sqrt{5})/2$ then for any odd prime $p$:
\begin{eqnarray*}
\mbox{if $p\equiv 1$ (mod $4$) then }&&\sum_{k=0}^{p-1}{2k\choose k}^2 {F_k\over 16^k}\equiv 0 \pmod{p^2} \\
\mbox{if $p\equiv 3$ (mod $4$) then }&&\sum_{k=0}^{p-1}{2k\choose k}^2 {L_k\over 16^k}\equiv 0 \pmod{p^2}. 
\end{eqnarray*}

\noindent (ii) The following result had been conjectured by Zhi-Wei Sun in \cite{Suzw:09}.

\noindent If $p\equiv 3$ (mod $4$),  by (ii) for $t=1/2$ we have
$$\sum_{k=0}^{p-1}{2k\choose k}^2 {1\over 32^k}\equiv 0 \pmod{p^2}$$
On the other hand, if $p\equiv 1$ (mod $4$) then let $f=(p-1)/4$ and let $x\equiv 1 \pmod{4}$ such that
$p=x^2+y^2$: Thus by our theorem  for $\alpha=1/2$ and $\beta=-1/2$ we obtain
$$\sum_{k=0}^{p-1}
{2k\choose k}^2 {1\over 32^k}\equiv
4^{-f}\sum_{k=0}^{2f}{2f\choose k}^2(-1)^k=
{(-4)^{-f}}{2f\choose f}\equiv 2x-{p\over 2x}\pmod{p^2}$$
where in the last step we used Theorem 9.4.3 in \cite{BEW:98}.
\section{Proof of the Theorem}
The following identity is well-known (see for example 3.18 in \cite{Go:72}):
for any non negative integer $n$, and any real numbers $\alpha,\beta$
$$\sum_{k=0}^{n}
{n\choose k}{n+k\choose k}(\alpha-\beta)^{n-k}\beta^k=
\sum_{k=0}^{n}{n\choose k}^2\alpha^{n-k}\beta^k.$$
We give here a straighforward proof via generating functions.
Let $F(z)=\sum_{k\geq 0}{n+k\choose k}z^k=(1-z)^{-(n+1)}$ then the LHS is a binomial transform and
therefore it is given by
$$[z^n]{1\over 1-(\alpha-\beta)z}\,F\left({\beta z\over 1-(\alpha-\beta)z)}\right)
={(1-(\alpha-\beta)z)^n\over (1-\alpha z)^{n+1}}.$$
Let $G(z)=\sum_{k\geq 0}{n\choose k}z^k=(1+z)^{n}$ then the RHS is another binomial transform and it is equal to
$$[z^n]{1\over 1-\alpha z}\,G\left({\beta z\over 1-\alpha z)}\right)
={(1-(\alpha-\beta)z)^n\over (1-\alpha z)^{n+1}}.$$
Thus the identity holds and since  
$${n\choose k}{n+k\choose k}={n+k\choose 2k}{2k\choose k}$$
it can be restated as follows
$$\sum_{k=0}^{n}
{n+k\choose 2k}{2k\choose k}(\alpha-\beta)^{n-k}\beta^k=
\sum_{k=0}^{n}{n\choose k}^2\alpha^{n-k}\beta^k.$$
In order to complete the proof it suffices to note that for $n=(p-1)/2$ 
we have that
$${n+k\choose 2k}={\prod_{j=1}^k(p^2-(2j-1)^2)\over 4^k(2k)!}\equiv
{\prod_{j=1}^k(2j-1)^2\over (-4)^k(2k)!}={2k\choose k}(-16)^{-k}\pmod{p^2}.$$
This key-congruence appears both in \cite{Suzw:09} (as an observation of Z. H. Sun) 
and in \cite{Mo:03} (see Proposition 2.5).
\hfill$\square$


\end{document}